\documentclass[11pt,reqno]{amsart}

\usepackage{color}
\usepackage{amssymb}
\usepackage{amsmath}
\usepackage{amsthm}
\usepackage{url}

\newtheorem{theorem}{Theorem}[section]
\newtheorem{definition}[theorem]{Definition}
\newtheorem{proposition}[theorem]{Proposition}

\newtheorem{conjecture}[theorem]{Conjecture}
\newtheorem{problem}[theorem]{Problem}

\begin{document}

\title[Derandomizing restricted isometries via the Legendre symbol]{Derandomizing restricted isometries\\via the Legendre symbol}

\author[Bandeira]{Afonso S.~Bandeira}
\address[Bandeira]{Program in Applied and Computational Mathematics, Princeton University, Princeton, NJ 08544, USA ({\tt ajsb@math.princeton.edu}).}

\author[Fickus]{Matthew Fickus}
\address[Fickus]{Department of Mathematics and Statistics, Air Force Institute of Technology, Wright-Patterson AFB, OH 45433, USA ({\tt matthew.fickus@afit.edu}).}

\author[Mixon]{Dustin G.~Mixon}
\address[Mixon]{Department of Mathematics and Statistics, Air Force Institute of Technology, Wright-Patterson AFB, OH 45433, USA ({\tt dustin.mixon@afit.edu}).}

\author[Moreira]{Joel Moreira}
\address[Moreira]{Department of Mathematics, Ohio State University, Columbus, OH 43210, USA ({\tt moreira.6@osu.edu}).}

\date{\today}

\begin{abstract}
The restricted isometry property (RIP) is an important matrix condition in compressed sensing, but the best matrix constructions to date use randomness.
This paper leverages pseudorandom properties of the Legendre symbol to reduce the number of random bits in an RIP matrix with Bernoulli entries.
In this regard, the Legendre symbol is not special---our main result naturally generalizes to any small-bias sample space.
We also conjecture that no random bits are necessary for our Legendre symbol--based construction.
\end{abstract}

\maketitle

\section{Introduction}

Erd\H{o}s first used the probabilistic method in 1947 to prove the existence of so-called \textit{Ramsey graphs}, i.e., $n$-vertex graphs with no clique or independent set of size $O(\log n)$~\cite{Erdos:47}.
Since then, the probabilistic method has been a popular technique for demonstrating the existence of various combinatorial objects, while at the same time leaving a desire for explicit constructions.
Even Ramsey graphs currently lack an explicit construction despite the fact that random graphs are (unverifiably) Ramsey with high probability.
This aggravatingly common state of affairs has been described by Avi Wigderson as the problem of \textit{finding hay in a haystack}.

Over the last decade, a new combinatorial object has been actively studied in applied mathematics:
An $M\times N$ matrix $\Phi$ is said to satisfy the \textit{$(K,\delta)$-restricted isometry property (RIP)} if
\begin{equation*}
(1-\delta)\|x\|_2^2
\leq\|\Phi x\|_2^2
\leq(1+\delta)\|x\|_2^2
\end{equation*}
for every vector $x$ with at most $K$ nonzero entries.
A matrix which satisfies RIP is well-suited as a \textit{compressive sensor}:
If $\Phi$ is $(2K,\delta)$-RIP with $\delta<\sqrt{2}-1$, then any nearly $K$-sparse signal $x$ can be stably reconstructed from noisy measurements $y=\Phi x+e$~\cite{Candes:08}.
Impressively, the number of measurements $M$ required for stable reconstruction is much smaller than the signal dimension~$N$.
As we will discuss in the next section, there are several known random constructions of RIP matrices with $M=O_\delta(K\operatorname{polylog}N)$, meaning the number of measurements can effectively scale with the complexity of the signal.

Unfortunately, RIP matrices provide yet another instance of finding hay in a haystack.
Much like Ramsey graphs, RIP matrices are ubiquitous, but verifying RIP is NP-hard in general~\cite{BandeiraDMS:13}.
To summarize the state of the art, almost every explicit construction of RIP matrices uses the following recipe (see~\cite{ApplebaumHSC:09,DeVore:07,FickusMT:12}, for example):
\begin{itemize}
\item[1.] Let the columns $\{\varphi_n\}_{n=1}^N$ of $\Phi$ have unit norm and satisfy $|\langle \varphi_n,\varphi_{n'}\rangle|\leq \mu$ for every $n\neq n'$.
\item[2.] Observe, either by the Riesz--Thorin theorem or the Gershgorin circle theorem, that every $K\times K$ principal submatrix $A$ of $\Phi^*\Phi$ satisfies $\|A-I\|_2\leq(K-1)\mu$.
\item[3.] Conclude that $\Phi$ is $(K,\delta)$-RIP for $\delta=(K-1)\mu$.
\end{itemize}
This coherence-based method of estimating $\delta$ requires $K\ll 1/\mu$.
However, Welch~\cite{Welch:74} gives the following universal bound 
\begin{equation*}
\mu^2\geq\frac{N-M}{M(N-1)},
\end{equation*}
which is $\geq1/(2M)$ whenever $M\leq N/2$.
As such, in Vinogradov notation we have
\begin{equation*}
M\gg \min\{1/\mu^2,N\}\gg\min\{K^2,N\},
\end{equation*}
meaning a coherence-based guarantee requires $M$ to scale with the \textit{square} of the complexity of the signal.
This inferiority to existing guarantees for random constructions compelled Terence Tao to pose the problem of explicitly constructing better RIP matrices on his blog~\cite{Tao:07}.
Since then, one explicit construction has managed to chip away at this power of two---Bourgain, Dilworth, Ford, Konyagin and Kutzarova~\cite{BourgainDFKK:11a,BourgainDFKK:11b} provided a new recipe, which is largely based on the following higher-order version of coherence:

\begin{definition}
We say a matrix $\Phi=[\varphi_1\cdots\varphi_N]$ has \textit{$(K,\theta)$-flat restricted orthogonality} if
\begin{equation*}
\bigg|\bigg\langle\sum_{i\in I}\varphi_i,\sum_{j\in J}\varphi_j\bigg\rangle\bigg|
\leq\theta(|I||J|)^{1/2}
\end{equation*}
for every disjoint pair of subsets $I,J\subseteq\{1,\ldots,N\}$ with $|I|,|J|\leq K$.
\end{definition}

\begin{theorem}[essentially Theorem~13 in~\cite{{BandeiraFMW:13}}, cf.~Lemma~1 in~\cite{BourgainDFKK:11a}]
\label{thm.fro to rip}
If a matrix $\Phi$ with unit-norm columns has $(K,\theta)$-flat restricted orthogonality, then $\Phi$ satisfies the $(2K,\delta)$-restricted isometry property with $\delta=150\theta\log K$.
\end{theorem}

Using this, Bourgain et al.\ leveraged additive combinatorics to construct an RIP matrix with $M=O_\delta(K^{2-\varepsilon})$ (according to~\cite{BourgainDFKK:11b}, $\varepsilon$ is on the order of $10^{-15}$); the reader is encouraged to read~\cite{Mixon:14} for a gentle introduction to this construction.

Since explicit constructions are evidently hard to come by, it is prudent to consider a relaxation of the hay-in-a-haystack problem:

\begin{problem}[Derandomized hay in a haystack]
Given $N$, $K$, $\delta$, and a budget $H$, what is the smallest $M$ for which there exists a random $M\times N$ matrix satisfying the $(K,\delta)$-restricted isometry property with high probability, but whose distribution only uses at most $H$ random bits?
\end{problem}

This derandomization problem is the focus of this paper.
The next section features a survey of the current literature.
In Section~3, we review our approach, namely, populating the matrix with consecutive Legendre symbols.
Our main result (Theorem~\ref{theorem_legendre}) gives that a random seed of $H=O_\delta(K\log K\log N)$ bits suffices to ensure that such a matrix satisfies $(2K,\delta)$-RIP with high probability; we also conjecture that no random seed is necessary (Conjecture~\ref{conj.main conjecture}).
Section~4 states a more general version of our main result---that any matrix with small-bias Bernoulli entries satisfies RIP with high probability (Theorem~\ref{thm.unbiased rip}); this is a natural generalization of a standard result with iid Bernoulli entries, and our main result follows from this one after verifying that random Legendre symbols have small bias.
We conclude in Section~5 with the proof of Theorem~\ref{thm.unbiased rip}.

\section{Prior work on derandomizing restricted isometries}

When $H=0$, the explicit construction of Bourgain et al.~\cite{BourgainDFKK:11a,BourgainDFKK:11b} holds the current record of $M=O_\delta(K^{2-\varepsilon})$, though this construction also requires $N\ll M^{1+\varepsilon}$.
If one is willing to spend $H=O_\delta(KN\log(N/K))$ random bits, then simply toss a fair coin for each entry of $\Phi$ and populate accordingly with $\pm1/\sqrt{M}$'s; this matrix is known to be $(K,\delta)$-RIP with high probability provided $M=O_\delta(K\log(N/K))$~\cite{BaraniukDDW:08,MendelsonPT:09}, and this is the smallest known $M$ for which there exists an $M\times N$ matrix satisfying $(K,\delta)$-RIP.
Far less entropy is required if one is willing to accept more structure in $\Phi$.
For example, toss $N$ independent weighted coins, each with success probability $M_0/N$, and then collect rows of the $N\times N$ discrete Fourier transform matrix which correspond to successful tosses; after scaling the columns to have unit norm, this random partial Fourier matrix will be $(K,\delta)$-RIP with high probability provided $M_0=O_\delta(K\log^3K\log N)$; see~\cite{CandesT:06,CheraghchiGV:13,RudelsonV:08}.
The entropy of this construction is $H=O(M_0\log(N/M_0))=O_\delta(K\log^3K\log^2N)$.

For another construction, this one requiring $N\leq M^2$, toss $M$ fair coins to populate the first column with $\pm1/\sqrt{M}$'s, and then take the remaining columns to be Gabor translates and modulates of this random seed vector; this matrix is $(K,\delta)$-RIP with high probability provided $M=O_\delta(K\log^2K\log^2N)$~\cite{KrahmerMR:14}, and uses entropy $H=M$.
In comparison with the random partial Fourier matrix, this construction has slightly more rows in exchange for a log factor of less entropy.
Another exchange can be made by modifying the partial Fourier construction:
Take a random partial Fourier matrix with $BM$ rows, partition the rows into subcollections of size $B$, toss $B$ fair coins for each subcollection and negate half of the rows accordingly, and then add the rows in each subcollection to produce a total of $M$ rows; after normalizing the columns, the resulting $M\times N$ matrix will be RIP with high probability provided $B=O(\log^{6.5} N)$ and $M=O_\delta(K\log^2K\log N)$~\cite{NelsonPW:14}.
Thanks to $B$, the entropy in this construction is $H=O_\delta(K\log^2K\log^{8.5}N)$, and so we save a log factor in $M$ in exchange for several log factors in $H$.

It should be noted that the preceding random matrix constructions were not developed with the specific intent of derandomization, but rather, they were developed to satisfy certain application-driven structural requirements (or to allow for fast matrix--vector multiplication) and derandomization came as a convenient byproduct.
Other constructions have also been proposed with the particular goal of derandomization.
For one example, start with a matrix of Bernoulli $0$'s and $1$'s, and then for each column, apply an error-correcting binary code to convert the column into a codeword and then change the $0$'s and $1$'s to $\pm1/\sqrt{M}$'s.
If the codebook covers $\{0,1\}^M$ in Hamming distance, then this matrix differs from the corresponding Bernoulli matrix of $\pm1/\sqrt{M}$'s in only a few entries, resulting in little loss in the $\delta$ for RIP.
This error-correction trick was introduced in~\cite{CalderbankJN:11}, which also identified various codes to derandomize Bernoulli matrices with $M=O_\delta(K\log(N/K))$.
Unfortunately, the total entropy for each of these cases is invariably $H=(M-\operatorname{polylog}M)N=\Omega(MN)$; the dual of concatenated simplex codes offers slightly more derandomization $H=(M-M^\varepsilon\log M)N=\Omega(MN)$, but at the price of $M=\Omega(K^{1+\varepsilon})$.

A particularly general RIP construction uses a \textit{Johnson--Lindenstrauss (JL) projection}, namely a random $M\times N$ matrix $\Phi$ such that for every $x\in\mathbb{R}^N$ of unit norm,
\begin{equation*}
\operatorname{Pr}\Big(\big|\|\Phi x\|_2^2-1\big|\geq \varepsilon\Big)
\leq e^{-cM\varepsilon^2}
\qquad
\forall \varepsilon>0.
\end{equation*}
JL projections are useful for dimensionality reduction, as they allow one to embed $P$ members of $\mathbb{R}^N$ in $M=O((\log P)/\varepsilon^2)$ dimensions while preserving pairwise distances to within a distortion factor of $\varepsilon$.
It was established in~\cite{BaraniukDDW:08} that JL projections are $(K,\delta)$-RIP with high probability provided $M=O_\delta(K\log(N/K))$.
Moreover, there are several random matrices which are known to be JL, for example, matrices of Bernoulli $\pm1/\sqrt{M}$'s~\cite{Achlioptas:03}.
Not surprisingly, there has also been work on derandomizing JL projections, though such projections inherently require $H=\Omega_\varepsilon(M)$; since $M<N$, then for every instance $\Phi_i$ of the random matrix $\Phi$, there exists a unit-norm vector $x$ in the null space of $\Phi_i$, and so the probability $p_i$ that $\Phi=\Phi_i$ is $\leq e^{-cM\varepsilon^2}$, implying
\begin{equation*}
H
=\sum_{i}p_i\log(1/p_i)
\geq\min_i\log(1/p_i)
\gg_\varepsilon M.
\end{equation*} 
This reveals an intrinsic bottleneck in derandomizing JL projections as a means to derandomize RIP matrices.
Nevertheless, the most derandomized RIP matrices to date are JL projections.

These particular JL projections exploit $r$-wise independent Bernoulli random variables, meaning every subcollection of size $r$ is mutually independent.
We note that the standard construction of $n$ Bernoulli random variables which are $r$-wise independent employs an $r$-wise independent hash family, requires only $O(r\log n)$ random bits, and is computationally efficient.
As a consequence of Theorem~2.2 in~\cite{ClarksonW:09}, a Bernoulli $\pm1/\sqrt{M}$ matrix with $\Omega_\varepsilon(M)$-wise independent entries is necessarily JL, thereby producing an RIP matrix with $M=O_\delta(K\log(N/K))$ and $H=O_\delta(K\log(N/K)\log N)$.
This has since been derandomized further by Kane and Nelson~\cite{KaneN:10}, who multiply a particular sequence of these JL projections to produce a JL family that leads to RIP matrices with $M=O_\delta(K\log(N/K))$ and $H=O_\delta(K\log(N/K)\log(K\log(N/K)))$; at the moment, this is the most derandomized construction with this minimal value of $M$.
In the spirit of derandomizing JL, note that Krahmer and Ward~\cite{KrahmerW:11} proved that randomly negating columns of an RIP matrix produces a JL projection, establishing that an RIP construction with entropy $H$ leads to a JL family with entropy $H+N$; of course, this is sub-optimal if $K$ is much smaller than $N$.

\section{Our approach}

In this paper, we propose a derandomized version of the matrix of Bernoulli $\pm1/\sqrt{M}$'s, and our construction makes use of the Legendre symbol.
For any integer $a$ and odd prime $p$, the \textit{Legendre symbol of $a$ and $p$} is defined to be
\begin{equation*}
\bigg(\frac{a}{p}\bigg)
:=
\left\{
\begin{array}{rl}
1&\mbox{if $\exists x\in\{1,\ldots,p-1\}$ such that $a\equiv x^2\bmod p$}\\
-1&\mbox{if $\nexists x\in\{0,\ldots,p-1\}$ such that $a\equiv x^2\bmod p$}\\
0&\mbox{if $a\equiv 0\bmod p$}.
\end{array}
\right.
\end{equation*}
The Legendre symbol is known for its pseudorandom behavior, and it is historically associated with derandomization.
For example, Erd\H{o}s used the probabilistic method in~\cite{Erdos:63} to show that almost every way of directing edges in a complete graph forms something called a \textit{paradoxical tournament}, and Graham and Spencer~\cite{GrahamS:71} later showed that directing edges between $p$ vertices according to the Legendre symbol (namely, $i\rightarrow j$ precisely when $((i-j)/p)=1)$ provides an explicit construction of such a tournament.
We take inspiration from this success in applying the Legendre symbol to find hay in a haystack.
Our particular application will make use of the following theorem (historically attributed to Harold Davenport), which establishes higher-order random-like cancellations:

\begin{theorem}[cf.\ Theorem~1 in~\cite{MauduitS:97}]
\label{thm.legendre_pseudorandom} 
There exists $p_0$ such that for every $p\geq p_0$ and integers $0<d_1<\cdots<d_k<p$ and $1\leq t\leq p-d_k$, we have
\begin{equation*}
\bigg|
\sum_{n=0}^{t-1}\bigg(\frac{n+d_1}{p}\bigg)\cdots\bigg(\frac{n+d_k}{p}\bigg)
\bigg|
\leq9kp^{1/2}\log p.
\end{equation*}
\end{theorem}

If the sum above ranged over all $n\in\mathbb{Z}/p\mathbb{Z}$, then we could invoke Weil's character sum estimate (Theorem 2C in~\cite{Schmidt:76}) to get $\leq kp^{1/2}$.
Instead, the sum is truncated so as to avoid a zero Legendre symbol, and so we are penalized with a log factor---this estimate is made possible by a version of the Erd\H{o}s--Tur\'{a}n inequality.

For one interpretation of this result, draw $X$ uniformly from $\{0,\ldots,p-k-1\}$.
Then all moments and mixed moments of the random Legendre symbols $\{((X+i)/p)\}_{i=1}^k$ have size $\leq (9kp^{1/2}\log p)/(p-k-1)$, which vanishes if $k=o(p^{1/2}/\log p)$.
As such, these random Legendre symbols have low correlation, but notice how little entropy we used.
Indeed, at the price of only $H<\log_2p$ random bits to pick $X$, we can enjoy as many as, say, $k=p^{1/3}$ bits of low correlation.
This quasi-paradox reveals something subtle about Bernoulli random variables:
While zero correlation implies independence (meaning the total entropy is the number of random variables), it is possible to produce many Bernoulli random variables of low correlation with very little entropy; as we discuss in the next section, this is related to packing many nearly orthogonal vectors in low-dimensional Euclidean space.

In this paper, we leverage the pseudorandomness identified in Theorem~\ref{thm.legendre_pseudorandom} to derandomize the matrix of Bernoulli $\pm1/\sqrt{M}$'s while maintaining RIP.
Our particular matrix construction uses the following recipe:
\begin{itemize}
\item[1.]
Given $N$, $K$ and $\delta$, pick $M$ and $H$ sufficiently large, and take $p$ to be some prime $\geq2^H+MN$.
\item[2.]
Draw $X$ uniformly from $\{0,\ldots,2^H-1\}$, and populate the entries of an $M\times N$ matrix $\Phi$ one column at a time with consecutive Legendre symbols $\{((X+i)/p)\}_{i=1}^{MN}$.
\item[3.]
Use Theorem~\ref{thm.legendre_pseudorandom} to show that $\Phi$ has flat restricted orthogonality, and therefore satisfies $(K,\delta)$-RIP by Theorem~\ref{thm.fro to rip}.
\end{itemize}
The reader should rest assured that the third step above is not obvious, but it summarizes our proof technique for the main result of this paper:

\begin{theorem}[Main Result]
\label{theorem_legendre}
Given $N$, $K$ and $\delta$, take
\begin{equation*}
M=(C_1/\delta^2)K\log^2 K\log N,
\qquad
H=C_2K\log((K/\delta)\log K)\log N
\end{equation*}
with $C_1$ and $C_2$ sufficiently large, and let $p$ denote any prime $\geq2^H+MN$.
Draw $X$ uniformly from $\{0,\ldots,2^H-1\}$, and define the corresponding $M\times N$ matrix $\Phi$ entrywise by the Legendre symbol
\begin{equation*}
\Phi[m,n]
:=\frac{1}{\sqrt{M}}\bigg(\frac{X+M(n-1)+m}{p}\bigg).
\end{equation*}
Then $\Phi$ satisfies the $(2K,\delta)$-restricted isometry property with high probability.
\end{theorem}

At this point, we explain how to efficiently construct this matrix.
Note that the output (i.e., the $M\times N$ matrix $\Phi$) is of length $MN\leq N^2$, and so consider a construction algorithm to be computationally efficient if it takes time polynomial in $N$.
The first task is to find a prime $p$ which is $\geq2^H+MN$.
By Bertrand's postulate, it suffices to search the integers between $2^H+MN$ and $2(2^H+MN)$, of which a fraction of about $1/\log(2^H+MN)=\Omega(1/H)$ are primes by the prime number theorem.
As such, we may randomly draw an integer in this interval, run the AKS primality test~\cite{AgrawalKS:04}, and repeat until we find a prime; this randomized algorithm will succeed in $\operatorname{poly}(H)=\operatorname{poly}(N)$ time with high probability.
(We note that deterministic alternatives to this approach have been studied, for example, in a polymath project~\cite{TaoCH:12}, but the best known algorithms of this form use superpolynomial time.)
Once the prime $p$ has been selected, we draw $X$ uniformly from $\{0,\ldots,2^H-1\}$, which requires $H=O(N)$ queries of a Bernoulli random number generator.
Finally, the entries of $\Phi$ are populated with $MN\leq N^2$ different Legendre symbols, and each of these can be calculated in $\operatorname{polylog}(p)=\operatorname{poly}(N)$ time, either by appealing to Euler's criterion (using repeated squaring with intermediate reductions modulo $p$) or by exploiting various identities of the Jacobi symbol (which equals the Legendre symbol whenever the ``denominator'' is prime)~\cite{Cohen:93}.

We now identify the shortcomings of our main result.
First, the number of measurements $M$ that scales like $K\log^2K\log N$ instead of the minimal $K\log(N/K)$.
We credit this (admittedly small) difference to our use of flat restricted orthogonality (i.e., Theorem~\ref{thm.fro to rip}) to demonstrate RIP; indeed, Theorem~14 in~\cite{{BandeiraFMW:13}} gives that using flat restricted orthogonality to prove RIP of Gaussian matrices, which are known to be RIP with minimal $M=O_\delta(K\log(N/K))$~\cite{CandesT:06}, leads to a scaling in $M$ that is identical to Theorem~\ref{theorem_legendre}, and so our result can perhaps be strengthened by somehow proving RIP ``directly.''
Next, we note that our construction uses entropy $H$ that scales like $K\log K\log N$, which is slightly more than the JL construction of Kane and Nelson~\cite{KaneN:10} that uses $O_\delta(K\log(N/K)\log(K\log(N/K)))$ random bits.
As we noted earlier, JL constructions of RIP matrices necessarily use at least $\Omega_\delta(M)=\Omega_\delta(K\log(N/K))$ random bits, whereas we believe our Legendre symbol construction can be derandomized quite a bit more:

\begin{conjecture}
\label{conj.main conjecture}
There exists a universal constant $C$ such that for every $\delta>0$, there exists $N_0>0$ and $P_0(N)=O(2^{\operatorname{poly}(N)})$ such that for every $(K,M,N)$ satisfying
\begin{equation*}
M\geq(C/\delta^2)K\log(N/K),
\qquad
N\geq N_0, 
\end{equation*}
and for every prime $p\geq P_0$, the $M\times N$ matrix $\Phi$ defined entrywise by the Legendre symbol
\begin{equation*}
\Phi[m,n]
:=\frac{1}{\sqrt{M}}\bigg(\frac{M(n-1)+m}{p}\bigg)
\end{equation*}
satisfies the $(2K,\delta)$-restricted isometry property.
\end{conjecture}

At its heart, Conjecture~\ref{conj.main conjecture} is a statement about how well the Legendre symbol exhibits additive cancellations.
In particular, if one is willing to use flat restricted orthogonality as a proof technique for RIP, the statement is essentially a bound on incomplete sums of Legendre symbols, much like those investigated in~\cite{Chung:94}.
We had difficulty proving these cancellations, and so we injected some randomness to enable the application of Theorem~\ref{thm.legendre_pseudorandom}.
Since Conjecture~\ref{conj.main conjecture} could very well be beyond the reach of modern techniques, we pose the following weaker and more general problem, whose solution would be a substantial advance in the derandomization of RIP (since this level of derandomization is not achievable with JL projections):

\begin{problem}[Breaking the Johnson--Lindenstrauss bottleneck]
\label{problem.JL bottleneck}
Find a construction of $M\times N$ matrices which satisfies the $(K,\delta)$-restricted isometry property with high probability whenever $M=O_{\delta}(K\operatorname{polylog}N)$ and uses only $H=o_{\delta;N\rightarrow\infty}(K\log(N/K))$ random bits.
\end{problem}

\section{The main result}

In this section, we reduce our main result to a more general (and independently interesting) result, which requires the following definition:

\begin{definition}
Let $\{X_i\}_{i=1}^n$ be a sequence of Bernoulli random variables taking values in $\{\pm1\}$.
We say $\{X_i\}_{i=1}^n$ is \textit{$\varepsilon$-biased} if 
\begin{equation}
\label{eq.epsilon biased}
\bigg|\mathbb{E}\prod_{i\in I}X_i\bigg|\leq\varepsilon
\end{equation}
for every nonempty $I\subseteq\{1,\ldots,n\}$.
\end{definition}

This notion of bias was introduced by Vazirani in~\cite{Vazirani:86}.
We will use bias as a means for derandomizing RIP matrices.
In particular, we will show that instead of using iid Bernoulli entries, it suffices to have small-bias Bernoulli entries for RIP, and then we will apply Theorem~\ref{thm.legendre_pseudorandom} to verify that random Legendre symbols have sufficiently small bias (even though the required entropy is small).
The idea that small-bias random variables can be used for derandomization is not new; indeed, one of their primary applications is to reduce the number of random bits needed for randomized algorithms~\cite{NaorN:90}.
As promised, the following generalizes the RIP result for matrices with independent Bernoulli entries to matrices with small-bias Bernoulli entries:

\begin{theorem}
\label{thm.unbiased rip}
Fix $C=5760000$.
There exists $N_0$ for which the following holds:
Given $N\geq N_0$, $K\geq1$ and $\delta>0$, take $M=(C/\delta^2)K\log^2 K\log N$ and pick
\begin{equation*}
\varepsilon
\leq\exp\Big(-40K\log((150/\delta)K\log K)\log N\Big).
\end{equation*}
Then any $M\times N$ matrix populated with $\varepsilon$-biased Bernoulli random variables (scaled by $1/\sqrt{M}$) satisfies the $(2K,\delta)$-restricted isometry property with probability $\geq1-2N^{-2K}$.
\end{theorem}

We note that the constants in this theorem have not been optimized.
Also, taking $\varepsilon=0$ recovers the iid Bernoulli result, though with extra log factors appearing in $M$.
We will prove this result in the next section.
For now, we establish that Theorem~\ref{thm.unbiased rip} implies our main result (Theorem~\ref{theorem_legendre}).
Drawing $X$ uniformly from $\{0,\ldots,2^H-1\}$, we bound the bias of the Legendre symbols $\{((X+i)/p)\}_{i=1}^{MN}$ using Theorem~\ref{thm.legendre_pseudorandom}:
\begin{equation*}
\bigg|\mathbb{E}\prod_{i\in I}\bigg(\frac{X+i}{p}\bigg)\bigg|
=\bigg|\frac{1}{2^H}\sum_{x=0}^{2^H-1}\prod_{i\in I}\bigg(\frac{x+i}{p}\bigg)\bigg|
\leq\frac{1}{2^H}|I|p^{1/2}\log p.
\end{equation*}
Since $|I|\leq MN$ and further $p\leq2(2^H+MN)\leq4\cdot2^H$ by Bertrand's postulate, we have
\begin{equation*}
\frac{1}{2^H}|I|p^{1/2}\log p
\leq\frac{1}{2^H}MN(4\cdot2^H)^{1/2}(\log 4+H\log 2)
\leq 4MNH\cdot 2^{-H/2}
\leq 4N^2\cdot2^{-H/3}.
\end{equation*}
As such, the random Legendre symbols we use to populate $\Phi$ are $(4N^2\cdot2^{-H/3})$-biased.
To use Theorem~\ref{thm.unbiased rip}, it remains to verify that this bias is sufficiently small:
\begin{equation*}
\log 4+2\log N-\frac{H}{3}\log 2\leq-40K\log((150/\delta)K\log K)\log N,
\end{equation*}
which is indeed satisfied by our choice of $H$ (i.e., taking $C_2$ to be sufficiently large).

Of course, we are not the first to use the Legendre symbol to produce small-bias random variables (for example, see~\cite{AlonGHP:92,Peralta:90}).
Also, there are several other constructions of small-bias random variables (e.g.,~\cite{AlonGHP:92,Ben-AroyaT:09,NaorN:90}), but these do not naturally lead to a conjecture such as Conjecture~\ref{conj.main conjecture}.
Interestingly, a stronger version of the Chowla conjecture~\cite{Chowla:65} implies that a randomly seeded portion of the Liouville function also produces small-bias random variables, and so one might pose a corresponding version of Conjecture~\ref{conj.main conjecture}.
Unfortunately, a large class of small-bias random variables cannot be used in conjunction with Theorem~\ref{thm.unbiased rip} to break the Johnson--Lindenstrauss bottleneck (Problem~\ref{problem.JL bottleneck}).
To see this, we first make an identification between small-bias random variables and linear codes.

A \textit{linear code} $\mathcal{C}\subseteq\mathbb{F}_2^q$ is a subspace of \textit{codewords}.
The $n\times q$ \textit{generator matrix} $G$ of an $n$-dimensional linear code $\mathcal{C}$ has the property that $\mathcal{C}=\{xG:x\in\mathbb{F}_2^n\}$, i.e., the rows of $G$ form a basis for $\mathcal{C}$.
The \textit{weight} of a codeword is the number of entries with value $1$.

\begin{proposition}[cf.~\cite{AlonGHP:92,Ben-AroyaT:09,NaorN:90}]
\
\begin{itemize}
\item[(a)]
Let $G$ be the generator matrix of an $n$-dimensional linear code in $\mathbb{F}_2^q$ such that every nonzero codeword has weight between $(1-\varepsilon)q/2$ and $(1+\varepsilon)q/2$.
Randomly sample $j$ uniformly over $\{1,\ldots,q\}$ and take $X_i:=(-1)^{G_{ij}}$ for every $i=1,\ldots,n$.
Then $\{X_i\}_{i=1}^n$ is $\varepsilon$-biased.
\item[(b)]
Let $\{X_i\}_{i=1}^n$ be an $\varepsilon$-biased sequence drawn uniformly from some multiset $\mathcal{X}\subseteq\{\pm1\}^n$ of size~$q$ ($\mathcal{X}$ is called an \textit{$\varepsilon$-biased set}).
For each $x\in\mathcal{X}$, consider the corresponding $g\in\mathbb{F}_2^n$ such that $x_i=(-1)^{g_i}$ for every $i=1,\ldots,n$.
Then the matrix $G$ whose columns are the $g$'s corresponding to $x$'s in $\mathcal{X}$ is the generator matrix of an $n$-dimensional linear code in $\mathbb{F}_2^q$ such that every nonzero codeword has weight between $(1-\varepsilon)q/2$ and $(1+\varepsilon)q/2$.
\end{itemize}
\end{proposition}

To see the significance of this proposition, consider the $n$-dimensional linear code $\mathcal{C}\subseteq\mathbb{F}_2^q$ corresponding to an $\varepsilon$-biased set of size $q$.
For each codeword $c\in\mathcal{C}$, define a unit vector $v_c\in\mathbb{R}^q$ whose entries have the form $(-1)^{c_i}/\sqrt{q}$.
It is easy to verify that $|\langle v_c,v_{c'}\rangle|\leq\varepsilon$ whenever $c\neq c'$, and so the Welch bound~\cite{Welch:74} gives
\begin{equation*}
\varepsilon^2
\geq\frac{2^n-q}{q(2^n-1)}.
\end{equation*}
For our application, we have $n=MN$ and $q=2^H$.
If $H\leq MN-1$ (i.e., $q\leq2^{n-1}$), then the Welch bound implies $\varepsilon^2\geq1/(2q)=2^{-(H+1)}$.
This gives the following result:

\begin{proposition}
Every $\varepsilon$-biased set $\mathcal{X}\subseteq\{\pm1\}^n$ has entropy $H\geq\min\{\log(1/\varepsilon),n-1\}$.
\end{proposition}

As such, applying Theorem~\ref{thm.unbiased rip} with any $\varepsilon$-biased set requires $H \gg_\delta K\log K\log N$.
In this sense, the Legendre construction in Theorem~\ref{theorem_legendre} is optimal.
On the other hand, no $\varepsilon$-biased set can be used with Theorem~\ref{thm.unbiased rip} to break the Johnson--Lindenstrauss bottleneck.

\section{Proof of Theorem~\ref{thm.unbiased rip}}

We will first show that $\Phi$ has $(K,\theta)$-flat restricted orthogonality, and then appeal to Theorem~\ref{thm.fro to rip} to get the $(2K,\delta)$-restricted isometry property.
To this end, fix a disjoint pair of subsets $I,J\subseteq\{1,\ldots,N\}$ with $|I|,|J|\leq K$.
We seek to bound the following probability:
\begin{equation}
\label{eq.prob to bound}
\operatorname{Pr}\left[\bigg|\bigg\langle\sum_{i\in I}\varphi_i,\sum_{j\in J}\varphi_j\bigg\rangle\bigg|>\theta(|I||J|)^{1/2}\right]
=\operatorname{Pr}\left[\bigg|\sum_{i\in I}\sum_{j\in J}\sum_{m=1}^M\Phi[m,i]\Phi[m,j]\bigg|>\theta(|I||J|)^{1/2}\right]
\end{equation}
Applying a version of Markov's inequality then gives
\begin{equation}
\label{eq.prob to bound 2}
\eqref{eq.prob to bound}
\leq\frac{1}{\big(\theta(|I||J|)^{1/2}\big)^q}\mathbb{E}\left[\bigg(\sum_{i\in I}\sum_{j\in J}\sum_{m=1}^M\Phi[m,i]\Phi[m,j]\bigg)^q\right],
\end{equation}
for some even integer $q$ which we will optimize later.
Observe that, because $q$ is even, the absolute value inside the expectation was not needed.
For now, we expand the product of sums and use linearity of expectation to get
\begin{equation*}
\eqref{eq.prob to bound 2}
=\frac{1}{\big(\theta(|I||J|)^{1/2}\big)^q}
\sum_{(i_1,j_1,m_1)\in I\times J\times [M]}\cdots\sum_{(i_q,j_q,m_q)\in I\times J\times [M]}
\mathbb{E}\bigg[\prod_{v=1}^q\Phi[m_v,i_v]\Phi[m_v,j_v]\bigg].
\end{equation*}
Next, since \eqref{eq.prob to bound 2} is nonnegative, the triangle inequality gives
\begin{equation*}
\eqref{eq.prob to bound 2}
=|\eqref{eq.prob to bound 2}|
\leq\frac{1}{\big(\theta(|I||J|)^{1/2}\big)^q}
\sum_{\{(i_v,j_v,m_v)\}_{v=1}^q\in(I\times J\times [M])^q}
\left|\mathbb{E}\bigg[\prod_{v=1}^q\Phi[m_v,i_v]\Phi[m_v,j_v]\bigg]\right|.
\end{equation*}

Most of the terms in the sum take the form of \eqref{eq.epsilon biased}, and so we can bound their contributions accordingly.
Each of the remaining terms has the property that all of its factors appear an even number of times in the product, and we will control the contribution of these terms by establishing how few they are.
With this in mind, we now define an index subset $S\subseteq(I\times J\times [M])^q$ (think ``surviving'' indices).
Explicitly, $\{(i_v,j_v,m_v)\}_{v=1}^q\in S$ if for every $v$, each of the sets
\begin{equation*}
\{v'\in[q]:(i_{v'},m_{v'})=(i_{v},m_{v})\},
\qquad
\{v'\in[q]:(j_{v'},m_{v'})=(j_{v},m_{v})\}
\end{equation*}
has an even number of elements.

We are able to use $\varepsilon$-bias to control all of the terms not in $S$.
Indeed, let $\{(i_v,j_v,m_v)\}_{v=1}^q\in S^\mathrm{c}$. Then there are $t\geq1$ distinct entries of $\Phi$ that appear in the product in (\ref{eq.prob to bound 2}) an odd number of times (observe that $I$ and $J$ are disjoint and hence an entry can not be simultaneously of the form $\Phi[m_v,i_v]$ and $\Phi[m_v,j_v]$).
Since there are a total of $2q$ entries of $\Phi$ in the product and the square of any entry is $1/M$, we can apply $\varepsilon$-bias to deduce that
\begin{equation}\label{eq_chowlatermsbound}
\left|\mathbb{E}\bigg[\prod_{v=1}^q\Phi[m_v,i_v]\Phi[m_v,j_v]\bigg]\right|
\leq M^{t/2-q}\frac{\varepsilon}{M^{t/2}}
=M^{-q}\varepsilon.
\end{equation}
Overall, the sum in our bound becomes
\begin{align*}
&\sum_{\{(i_v,j_v,m_v)\}_{v=1}^q\in S}\left|\mathbb{E}\bigg[\prod_{v=1}^q\Phi[m_v,i_v]\Phi[m_v,j_v]\bigg]\right|
+\sum_{\{(i_v,j_v,m_v)\}_{v=1}^q\in S^\mathrm{c}}\left|\mathbb{E}\bigg[\prod_{v=1}^q\Phi[m_v,i_v]\Phi[m_v,j_v]\bigg]\right|\\
&\qquad\qquad \leq M^{-q}\left(|S|+\varepsilon|S^\mathrm{c}|\right),
\end{align*}
where the inequality applies the fact that $|\Phi[m,n]|=1/\sqrt{M}$ to each term from $S$ and (\ref{eq_chowlatermsbound}) to each term from $S^\mathrm{c}$.
At this point, we wish to bound $|S|$.

For each $s=\{(i_v,j_v,m_v)\}_{v=1}^q\in S$, there exist two perfect matchings of $\{1,\ldots,q\}$, say $\mathcal{M}_1$ and $\mathcal{M}_2$, such that
\begin{equation*}
(i_{v},m_{v})=(i_{v'},m_{v'})~\forall\{v,v'\}\in \mathcal{M}_1,
\quad
(j_{v},m_{v})=(j_{v'},m_{v'})~\forall\{v,v'\}\in \mathcal{M}_2.
\tag{$P[s,\mathcal{M}_1,\mathcal{M}_2]$}
\end{equation*}
(Here, we use the phrase ``perfect matching'' so as to convey that $\{1,\ldots,q\}$ is partitioned into sets of size $2$.
Let $\mathcal{M}(q)$ denote all such perfect matchings.)
We now start our bound:
\begin{align*}
|S|
&=\#\Big\{s\in(I\times J\times [M])^q:\exists\mathcal{M}_1,\mathcal{M}_2\in\mathcal{M}(q)\mbox{ such that }P[s,\mathcal{M}_1,\mathcal{M}_2]\Big\}\\
&\leq\#\Big\{(s,\mathcal{M}_1,\mathcal{M}_2):P[s,\mathcal{M}_1,\mathcal{M}_2]\Big\}\\
&=\sum_{\mathcal{M}_1\in\mathcal{M}(q)}\sum_{\mathcal{M}_2\in\mathcal{M}(q)}\#\Big\{s:P[s,\mathcal{M}_1,\mathcal{M}_2]\Big\}.
\end{align*}
Next, given a perfect matching $\mathcal{M}\in\mathcal{M}(q)$, let $C(\mathcal{M})$ denote the set of $q$-tuples $\{m_v\}_{v=1}^q$ satisfying $m_{v}=m_{v'}$ whenever $\{v,v'\}\in \mathcal{M}$ (these are the $q$-tuples which consistently ``color'' the matching).
Note that the $q$-tuple of $m_v$'s from $s=\{(i_v,j_v,m_v)\}_{v=1}^q$ satisfying $P[s,\mathcal{M}_1,\mathcal{M}_2]$ lies in both $C(\mathcal{M}_1)$ and $C(\mathcal{M}_2)$.
Moreover, the $q$-tuple of $i_v$'s from $s=\{(i_v,j_v,m_v)\}_{v=1}^q$ satisfying $P[s,\mathcal{M}_1,\mathcal{M}_2]$ assigns members of $I$ to the $q/2$ pairs $\{v,v'\}\in\mathcal{M}_1$, meaning the $q$-tuple is one of $|I|^{q/2}$ such possibilities; similarly, there are $|J|^{q/2}$ possibilities for the $q$-tuple of $j_v$'s.
As such, we can continue our bound:
\begin{equation*}
|S|
\leq|I|^{q/2}|J|^{q/2}\sum_{\mathcal{M}_1\in\mathcal{M}(q)}\sum_{\mathcal{M}_2\in\mathcal{M}(q)}|C(\mathcal{M}_1)\cap C(\mathcal{M}_2)|.
\end{equation*}
At this point, we note that $\sum_{\mathcal{M}_2\in\mathcal{M}(q)}|C(\mathcal{M}_1)\cap C(\mathcal{M}_2)|$ does not depend on $\mathcal{M}_1$, since any permutation of $\{1,\ldots,q\}$ that sends one version of $\mathcal{M}_1$ to another can be viewed as a permutation acting on $\mathcal{M}_2$ (thereby merely permuting the summands).
As such, we can express our bound in terms of a ``canonical matching'' $\mathcal{M}_0:=\{\{1,2\},\{3,4\},\ldots,\{q-1,q\}\}$:
\begin{equation*}
|S|
\leq|I|^{q/2}|J|^{q/2}|\mathcal{M}(q)|\sum_{\mathcal{M}_2\in\mathcal{M}(q)}|C(\mathcal{M}_0)\cap C(\mathcal{M}_2)|.
\end{equation*}
Note that $|\mathcal{M}(q)|=(q-1)!!$ since there are $q-1$ possible members to match with $1$, and then $q-3$ members remaining to match with the smallest free element, etc.
It remains to determine the sum over $\mathcal{M}(q)$.
To parse this sum, consider a graph of $q$ vertices with edge set $\mathcal{M}_0\cup\mathcal{M}_2$.
Then $|C(\mathcal{M}_0)\cap C(\mathcal{M}_2)|$ counts the number of ways of coloring the vertices using $M$ colors in such a way that each component (which is necessarily a cycle) has vertices of a common color.
We will use this interpretation to prove the following claim:
\begin{equation*}
\sum_{\mathcal{M}_2\in\mathcal{M}(q)}|C(\mathcal{M}_0)\cap C(\mathcal{M}_2)|
=\frac{(M+q-2)!!}{(M-2)!!}.
\end{equation*}
To proceed, start with a graph $G$ of $q$ vertices, and let the edge set be $\mathcal{M}_0$.
Let $A$ denote the ``available'' set of vertices, and initialize $A$ as $\{1,\ldots,q\}$.
We will iteratively update the pair $(G,A)$ as we build all possible matchings $\mathcal{M}_2$ and colorings.
At each step, denote $k:=|A|$.
Take the least-numbered available vertex and draw an edge between it and any of the other $k-1$ available vertices.
Add this edge to $\mathcal{M}_2$ and remove both vertices from $A$.
If this edge completes a cycle in $G$, then color that cycle with any of the $M$ colors.
Overall, in this step you either pick one of $k-2$ vertices which do not complete the cycle, or you complete the cycle and pick one of $M$ colors; this totals to $M+k-2$ choices at each step.
Since the number $k$ of available vertices decreases by $2$ at each step, the product of these totals corresponds to the claim.

At this point, we summarize the status of our bound:
\begin{align}
\nonumber
\operatorname{Pr}\left[\bigg|\bigg\langle\sum_{i\in I}\varphi_i,\sum_{j\in J}\varphi_j\bigg\rangle\bigg|>\theta(|I||J|)^{1/2}\right]
&\leq\frac{1}{\big(M\theta(|I||J|)^{1/2}\big)^q}\bigg(|S|+\varepsilon|S^\mathrm{c}|\bigg)\\
\label{eq.two terms to bound}
&\leq\frac{(q-1)!!(M+q-2)!!}{(M\theta)^q(M-2)!!}+\frac{\varepsilon|I|^{q/2}|J|^{q/2}}{\theta^q},
\end{align}
where the last inequality applies our bound on $|S|$ along with $|S^\mathrm{c}|\leq|(I\times J\times [M])^q|=(|I||J|M)^q$.
To continue, we will bound the terms in \eqref{eq.two terms to bound} separately.
When $M$ is even (otherwise, $M-1$ is even and similar bounds will apply), the first term has the following bound:
\begin{align*}
\frac{(q-1)!!(M+q-2)!!}{(M\theta)^q(M-2)!!}
&\leq\frac{q!!(M+q-2)!!}{(M\theta)^q(M-2)!!}\\
&\leq\frac{2^q((\frac{q}{2})!)^2}{(M\theta)^q}\binom{\frac{M}{2}+\frac{q}{2}-1}{\frac{q}{2}}
\leq\frac{2^qe^2(\frac{q/2+1}{e})^{q+2}}{(M\theta)^q}\bigg(e\cdot\frac{\frac{M}{2}+\frac{q}{2}-1}{\frac{q}{2}}\bigg)^{q/2},
\end{align*}
where the last step uses Stirling-type bounds on both the factorial and the binomial coefficient.
Next, we combine like exponents to get
\begin{align*}
\frac{(q-1)!!(M+q-2)!!}{(M\theta)^q(M-2)!!}
&\leq e^2\Big(\frac{q/2+1}{e}\Big)^2\bigg(\frac{4}{(M\theta)^2}\Big(\frac{q/2+1}{e}\Big)^2\cdot e\cdot\Big(\frac{M}{q}+1-\frac{2}{q}\Big)\bigg)^{q/2}\\
&\leq q^2\bigg(\frac{4}{(M\theta)^2}\Big(\frac{q}{e}\Big)^2\cdot e\cdot\frac{2M}{q}\bigg)^{q/2},
\end{align*}
where the last step assumes $q\geq 2$ (given our choice of $q$ later, this will occur for sufficiently large $N$) so that $q/2+1\leq q$ and also assumes $q\leq M$ (which will follow from the fact that $\delta\leq 1$) so that $M/q+1\leq 2M/q$.
At this point, we take $q=M\theta^2/8$ (or more precisely, take $q$ to be the smallest even integer larger than $M\theta^2/8$) to get
\begin{equation*}
\frac{(q-1)!!(M+q-2)!!}{(M\theta)^q(M-2)!!}
\leq \Big(\frac{M\theta^2}{8}\Big)^2\cdot e^{-M\theta^2/16}.
\end{equation*}
This can be simplified by noting that $M=(C/\delta^2)K\log^2K\log N$ and $\theta=\delta/(150\log K)$ together give $M\theta^2=256K\log N$:
\begin{equation}
\label{eq.first term bound}
\frac{(q-1)!!(M+q-2)!!}{(M\theta)^q(M-2)!!}
\leq \operatorname{exp}\bigg(2\log\Big(\frac{M\theta^2}{8}\Big)-\frac{M\theta^2}{16}\bigg)
\leq\operatorname{exp}\Big(-8K\log N\Big),
\end{equation}
for sufficiently large $N$.
Next, we bound the second term in \eqref{eq.two terms to bound}:
\begin{equation}
\label{eq.second term bound}
\frac{\varepsilon|I|^{q/2}|J|^{q/2}}{\theta^q}
\leq\varepsilon(K/\theta)^q
=\varepsilon((150/\delta)K\log K)^q
\leq\exp\Big(-8K\log((150/\delta)K\log K)\log N\Big)
\end{equation}
Now we can combine the bounds \eqref{eq.first term bound} and \eqref{eq.second term bound} to get
\begin{align*}
&\operatorname{Pr}\left[\bigg|\bigg\langle\sum_{i\in I}\varphi_i,\sum_{j\in J}\varphi_j\bigg\rangle\bigg|>\theta(|I||J|)^{1/2}\right]\\
&\qquad\leq \operatorname{exp}\Big(-8K\log N\Big)+\exp\Big(-8K\log((150/\delta)K\log K)\log N\Big)\\
&\qquad\leq 2\operatorname{exp}\Big(-8K\log N\Big).
\end{align*}

At this point, we are ready to show that $\Phi$ has the $(2K,\delta)$-restricted isometry property with high probability.
This calculation will use a union bound over all possible choices of disjoint $I,J\subseteq\{1,\ldots,N\}$ with $|I|,|J|\leq K$, of which there are
\begin{equation*}
\sum_{|I|=1}^K\sum_{|J|=1}^K\binom{N}{|I|}\binom{N-|I|}{|J|}
\leq K^2\binom{N}{K}^2
\leq K^2\Big(\frac{eN}{K}\Big)^{2K}
\leq\operatorname{exp}\Big(4K+2K\log(N/K)\Big).
\end{equation*}
With this, we apply a union bound to get
\begin{align*}
\operatorname{Pr}\Big[\mbox{$\Phi$ is not $(2K,\delta)$-RIP}\Big]
&\leq\operatorname{Pr}\Big[\mbox{$\Phi$ does not have $(K,\theta)$-FRO}\Big]\\
&\leq\operatorname{exp}\Big(4K+2K\log(N/K)\Big)\cdot\operatorname{Pr}\left[\bigg|\bigg\langle\sum_{i\in I}\varphi_i,\sum_{j\in J}\varphi_j\bigg\rangle\bigg|>\theta(|I||J|)^{1/2}\right]\\
&\leq2\operatorname{exp}\Big(4K+2K\log(N/K)-8K\log N\Big)\\
&\leq2\operatorname{exp}\Big(-2K\log N\Big),
\end{align*}
which completes the result.

\section*{Acknowledgments}

The authors thank Prof.\ Peter Sarnak for insightful discussions.
A.\ S.\ Bandeira was supported by AFOSR Grant No.\ FA9550-12-1-0317, and M.\ Fickus and D.\ G.\ Mixon were supported by NSF Grant No.\ DMS-1321779.
The views expressed in this article are those of the authors and do not reflect the official policy or position of the United States Air Force, Department of Defense, or the U.S. Government.

\end{document}